\documentclass{aastex}          
\usepackage{spr-astr-addons}    

\begin{document}
%
\title{On a planar circular restricted charged three-body problem}

\shorttitle{Restricted charged three-body problem}
\shortauthors{A. Bengochea, C. Vidal}


\author{Abimael Bengochea} \affil{Departamento de Matem\'aticas, UAM Iztapalapa. Av. San Rafael Atlixco 186, M\'exico, D.F. 09340} \email{abc@xanum.uam.mx}

\author{Claudio Vidal} \affil{Departamento de Matem\'atica, Facultad de
Ciencias, Universidad  del Bio Bio. Casilla 5-C, Concepci\'on,
VIII-Regi\'on, Chile} \email{clvidal@ubiobio.cl}


\begin{abstract}
We introduce a circular restricted charged three-body problem on the plane. In this model, the gravitational 
and Coulomb forces, due to the primary bodies, act on a test particle; the net force exerted by some primary body 
on the test particle can be attractive, repulsive or null. The restricted problem is
obtained by the general planar charged three-body problem considering one mass of the three bodies going to zero. We 
obtain necessary restrictions for the parameters that appear in the problem, in order to be
well defined. Taking into account such restriction, we study the existence and linear
stability of the triangular equilibrium solutions, as well as its
location in the configuration space. We also obtain necessary and sufficient conditions for the existence of the collinear equilibrium solutions.
\end{abstract}

\keywords{Hamiltonian Vector Fields, charged three-body problem,
circular restricted charged three-body problem.}

\section{Introduction}

The force exerted by electric ${\bf  E}$ and magnetic field
${\bf B}$ on a particle of velocity ${\bf v}$ and charge $q$ is
${\bf F} = q ({\bf E} + {\bf v} \times {\bf B})$, and it is denominated
Lorentz's force \citep{goldstein, jackson}. The exact description of this problem is formulated in the
relativistic context. Nevertheless, when the speed of the particle
is much smaller than the one of the light $c$ (around
$300,000km/s$), the dominant term is the Coulombian since the
first relativistic correction is of order $ (v/c)^2$, as it is
shown by the Lagrangian of Darwin \citep{jackson}. In the study of charged particles, the gravitational
interaction is commonly of smaller magnitude than the electromagnetic one, whereas in celestial
problems it happens the opposite.

Several works have been written about restricted models concerning three point particles (it is assumed that one of the 
bodies does not affect the motion of the other two particles) with additional interactions to the gravitational 
one, for instance the Coulomb or photogravitational cases (the gravitational force is attractive, the Coulomb can be either attractive or 
repulsive, and the photogravitational only repulsive, in a generic way - the zero value is allowed in all of them). 

\cite{rad-1,rad-2} introduced the photogravitational restricted three-body problem, and \cite{dionysiou-2} a restricted charged three-body problem. \cite{schuerman} studied the stability, and location on the configuration space, of the equilibrium solutions when radiation pressure and Poynting-Roberts forces are included. In the photogravitational case, \cite{kunitsyn-1} studied the existence and stability of the collinear equilibrium solutions, \cite{lukyanov-1} the existence of the triangular equilibrium solutions, and the collinear equilibrium solution, in a parametric way. In the same model, \cite{kunitsyn-2} considered the existence of the collinear equilibrium solutions, and \cite{kunitsyn-3} the existence and stability of the equilibrium triangular solutions. \cite{simmons} studied the restricted case with radiation pressure; part of his study is concerned with the analysis of the triangular, collinear and spatial equilibrium solutions. In the photogravitational restricted problem, \cite{lukyanov-2} studied the stability of the collinear and triangular equilibrium solutions, as well as its location in the configuration space.

Some of the studies present the same potential function, including our own. However, the parameters, as well as its allowed values, can be different. For instance, in the restricted charged three-body problem, \cite{dionysiou-1} consider a different parameterization from
the one used in this work. In fact, the potential function is
$V=\frac{q-\mu}{\rho_1}+
 \frac{\mu-q}{\rho_2}$, where $\mu \in (0,1/2]$, $q \in \mathbb{R}$, and $\rho_i$ denotes
the distance between the body $i=1,2$ and the third body. On the other hand, in the photogravitational case, \cite{rad-1}, \cite{kunitsyn-1},
\cite{lukyanov-1} and \cite{simmons}, consider a parameterization similar to ours. The potential function is
$V=\frac{\gamma_1(1-\mu)}{\rho_1}+\frac{\gamma_2
\mu}{\rho_2}$, where $\mu$, $\gamma_1$ and $\gamma_2$ are real
parameters, and $\rho_i$, $i=1,2$, the distance between the bodies, as in the charged case. The parameter $\mu$ belongs to $(0,1/2]$, whereas $\gamma_i \in (-\infty, 1]$, $i=1,2$. The qualitative properties of the gravitational 
and photogravitational forces acting on the test body, by the primaries $i = 1,2$, are defined by the 
parameter $\gamma_i$: if $\gamma_i<0$ the gravitational force dominates over the 
photogravitational one, for $\gamma_i=0$ both forces are equal to each other in 
magnitude, for $0<\gamma_i<1$ the gravitational is the strongest one, and for $\gamma_i=1$ the 
photogravitational force is zero. The case $\gamma_i>1$ has not physical sense because it 
corresponds to an attractive photogravitational force; usually it is included
in the study of the mathematical model. The
different selection of
units and generality with which these problems have been studied
can be appreciated in the number of parameters that appear, that is, two and three respectively.

This paper is organized as follows. In Section 2 we give a brief review
of the problem of two charged bodies. Here we emphasize that the necessary
condition on the parameters of mass $m_i$ and charge $q_i$, $i=1,2$, for the existence of circular orbits of the
particles $1$ and $2$, is $G- k \ \frac{q_1q_2}{m_1m_2} > 0$ ($G$ is the
constant of universal gravitation and $k$ is the constant of
Coulomb); this limits the values of the parameters that appear in the restricted problem of three charged bodies on 
the plane, which we enunciate in Section 3. With this aim, we make the
mass of the third body tends to zero, which gives us a well
defined restricted circular problem. The
potential function associated to this problem is given by
$V = \frac{\beta_1 (1- \mu)}{\rho_1} + \frac{\beta_2 \
\mu}{\rho_2} $, with parameters $\mu$, $\beta_1$, $\beta_2$, and the usual distances
$\rho_1$, $\rho_2$. Here $\mu \in (0,1/2]$ and $\beta_i \in \mathbb{R}$ for $i=1,2$, with the
restriction $(\beta_1-1) (\beta_2-1) < 1$, which is the necessary condition for having circular solutions
for the bodies $1$ and $2$. Such condition has not been considered in previous studies, for instance \cite{dionysiou-1},
\cite{dionysiou-2}, \cite{kunitsyn-1}, \cite{kunitsyn-3}, \cite{lukyanov-1}, \cite{lukyanov-2},
\cite{simmons}, reason why our problem acquires a great difference
to others already treated and justifies its study. We remark that, if the parameters that
appear in the potential are not properly considered, the interaction between the primary bodies
will be repulsive and might prevent those bodies from having a circular movement. 
As in the photogravitational model, the parameters $\beta_i$, $i = 1,2$ determine the qualitative features of the gravitational and Coulomb forces. Actually, we have a similar description for both photogravitational and charged problems, with the exception that $\beta_i > 1$ is physically possible and corresponds to a Coulomb attractive force. Once
established the restricted charged three-body problem, in Section 4 we study
the existence of triangular and collinear equilibrium solutions. In Section 5 we deal with the linear stability of the triangular equilibrium solutions and its location in the configuration and parameters space. We finish the article with the conclusions of this work.

\section{Dynamics of the two-charged problem}

Consider an inertial frame of reference. The Hamiltonian associated to bodies $1$ and $2$ of charges $q_1$ and $q_2$, and
masses $m_1$ and $m_2$ respectively, with gravitational and Coulombian interaction, is
\begin{equation}
H({\bf r}_1, {\bf r}_2, {\bf p}_1, {\bf p}_2)=\frac{\|{\bf
p}_1\|^2}{2m_1} + \frac{\|{\bf p}_2\|^2}{2m_2} +
 \frac{kq_1q_2-Gm_1m_2}{\|{\bf r}_1 - {\bf r}_2 \|},
 \label{eq-1}
\end{equation}
where $k > 0$ is the Coulomb constant, $G > 0$ is the universal
gravitational constant, and ${\bf r}_i \in \mathbb{R}^3$, ${\bf p}_i \in \mathbb{R}^3$, $i=1,2$ are the positions and moments of the bodies with the same index, respectively.

\begin{figure}[!h]
\centering
\includegraphics[scale=0.6]{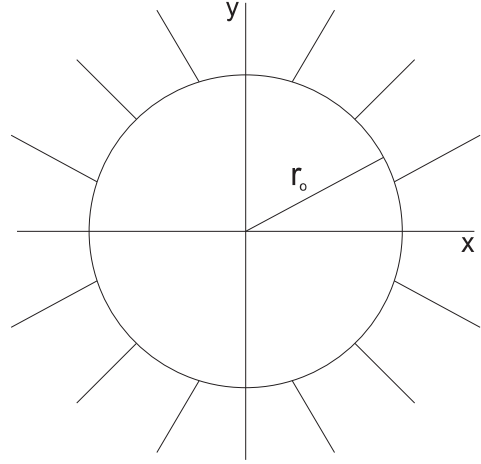}
\caption{Hill's region (hatched part) of the two-charged problem, with $h_{rel}^*
> 0$, $C < 0$.}
\label{claudio_fig1}
\end{figure}

Given the Hamiltonian (\ref{eq-1}), and the initial
conditions, the equations of Hamilton $\dot{{\bf
r}_i} =H_{{\bf p}_i}, \dot{{\bf p}_i} = - H_{{\bf r}_i}$ determine the motion
of each body. The equations of motion correspond to
a system of differential equations of second order:
$$
 \begin{array}{l}
\displaystyle m_1 \ddot{\bf r}_1=  -(G m_1 m_2 - k q_1 q_2) \frac{{\bf r}_1-
 {\bf r}_2}{\|{\bf r}_1-
 {\bf r}_2\|^3}, \\[1pc]
\displaystyle m_2 \ddot{\bf r}_2=  -(G m_1 m_2 - k q_1 q_2) \frac{{\bf r}_2-
 {\bf r}_1}{\|{\bf r}_1-
 {\bf r}_2\|^3},
 \label{sis-2-cargas-1}
 \end{array}
$$
where the dot denotes derivative with respect to time $t$. In order to classify the solutions, we define 
$C=G m_1 m_2-k q _1 q_2$. Three cases are identified:
\begin{itemize}
\item $C>0$,

\item $C=0$,

\item $C<0$.
\end{itemize}
The first one is equivalent to the Kepler's problem whose
dynamics is well known. Notice that the condition $C
> 0$ is satisfied whenever the charges have different signs, and
if the charges have equal signs it is required that $k q_1q_2 <
Gm_1m_2$. The case $C = 0$ corresponds to two free particles. The case $C<0$ is rarely discussed in the 
literature, therefore we give a brief description of it, following the same steps as in the Kepler's problem \citep{goldstein}. Notice that the system associated to $C <0$ has the Kepler's constants of motion, namely energy, angular momentum and those related to the center of mass of the system. Therefore, in a generic way, the movement of the two bodies occurs in a fixed plane, which we assume from here on. As first step, we make a symplectic transformation of coordinates, that change the
vectors ${\bf r}_1$, ${\bf r}_2$, ${\bf p}_1$, ${\bf p}_2$ to the ones
$$
\begin{array}{c}
\displaystyle {\bf r} = {\bf r}_2- {\bf r}_1, \quad  {\bf R} = \frac
{m_1 {\bf r}_1+ m_2 {\bf r}_2} {m_1+m_2}, \\ [1pc]
\displaystyle {\bf p} = \frac{m_1}{m_1 + m_2} {\bf p}_2 - \frac{m_2}{m_1 + m_2} {\bf p}_1, \quad {\bf P} = {\bf p}_1 + {\bf p}_2.
\end{array}
$$
The upper vectors correspond to the center of mass of the system, whereas the lowercase ones describe the motion relative to the center of mass of the system. From the previous change of variables it is obtained
$$
\begin{array}{l}
\displaystyle {\bf r}_1=  -\frac{m_2}{m_1+m_2} {\bf r} + {\bf R}, \quad {\bf r}_2=  \frac{m_1}{m_1+m_2} {\bf r} + {\bf R}, \\[1pc]
\displaystyle {\bf p}_1= -{\bf p} + \frac{m_1}{m_1+m_2} {\bf P}, \quad {\bf p}_2= {\bf p} + \frac{m_2}{m_1+m_2} {\bf P}. \label{ecuacion-2}
\end{array}
$$
By replacing the new variables and $-C = \vert C \vert$ in the
Hamiltonian (\ref{eq-1}) we obtain
\begin{equation}
H^*({\bf r}, {\bf R}, {\bf p}, {\bf P})=\frac{\|{\bf p}\|^2}{2\mu}
+ \frac{\|{\bf P}\|^2}{2M} + \frac{\vert C \vert}{\|{\bf r\|}},
\label{eq-2}
\end{equation}
where $\mu = \frac{m_1m_2}{m_1+m_2}$ is the reduced mass, and
$M=m_1+m_2$ is the total mass of the system. The Hamiltonian (\ref{eq-2}) is separable. For the motion relative to the center of mass we have
\begin{equation}
H_{rel}^*({\bf r}, {\bf p})=\frac{\|{\bf p}\|^2}{2\mu} +
\frac{\vert C \vert}{\|{\bf r}\|}.
\label{eq-3}
\end{equation}
It is not difficult to identify the Hill's region for (\ref{eq-3}). This consists
of all the points ${\bf r} $ on the configuration space such that
$$
 h_{rel}^* \ge \frac{\vert C \vert}{\|{\bf r}\|},
\label{eq-4}
$$
where $h_{rel}^*>0$ is the constant of motion associated
to (\ref{eq-3}). Thus, the points that
define this region must satisfy
$$
\|{\bf r}\| \geq r_0,
\label{eq-5}
$$
with $r_0 = \frac{\vert C \vert}{h_{rel}^*}$. The Hill's region is outlined in Figure~\ref{claudio_fig1}.

\placefigure{claudio_fig1}

\begin{figure}[!h]
\centering
\includegraphics[scale=0.8]{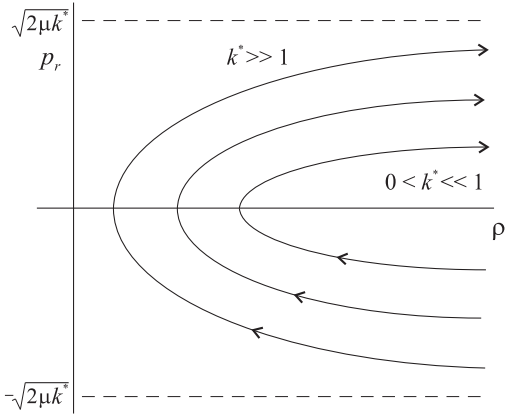}
\caption{Phase portrait of the two-charged problem, with $k^*
> 0, C < 0$.}
\label{claudio_fig2}
\end{figure}

Now we introduce polar coordinates by means of using the generating function
$S_2(p_x,p_y,\rho,\theta)=p_x \rho \cos (\theta) + p_y \rho \sin
(\theta)$:
$$
\begin{array}{c}
x = \frac{\partial S_2}{\partial p_x} = \rho \cos (\theta), \\[1pc]
y = \frac{\partial S_2}{\partial p_y} = \rho \sin (\theta), \\[1pc]
p_r = \frac{\partial S_2}{\partial \rho} = p_x \cos (\theta) + p_y \sin (\theta), \\[1pc]
p_{\theta} = \frac{\partial S_2}{\partial \theta} = -p_x \rho \sin (\theta) + p_y \rho \cos (\theta).
\end{array}
$$
With this transformation we obtain the new Hamiltonian
\begin{equation}
K(\rho,\theta,p_r,p_\theta)=\frac{1}{2\mu}\left(p_r^2+\frac{p_\theta^2}{\rho^2}\right)
+ \frac{\vert C \vert}{\rho},
\label{eq-4}
\end{equation}
associated to the constant of movement $k^*=h_{rel}^*>0$. In this Hamiltonian does not appear $\theta$, reason why the
canonical momentum $p_ {\theta} =l$ is constant, which corresponds
to the conservation of the angular momentum. Replacing $p_{\theta}
=l$ in (\ref{eq-4}), and considering the reduced energy $k^*$, it is obtained
\begin{equation}
k^*= \frac{p_r^2} {2 \mu} + \frac{l^2} {2 \mu \rho^2} + \frac
{\vert C \vert} {\rho},
\label{eq-5}
\end{equation}
where $\frac{l^2} {2 \mu \rho^2} + \frac {\vert C \vert}{\rho}$ is defined
as the effective potential $V_{\it eff}(\rho)$. The constant of movement (\ref{eq-5}) is useful to obtain
the phase portrait; from this we get $p_r= \pm
\sqrt{2 \mu (k^* - V_{\it eff}(\rho))}$. When $\rho$ takes the value
$\rho_ {*} = \frac{\vert C \vert}{2k^*} (1+ \sqrt{1+
\frac{2l^2k^*}{C^2 \mu}})$ the canonical momentum $p_r$ is zero.
Then, as $k^*$ grows $\rho_{*} $ tends to
zero, and if $\rho \to \infty$ then $p_r \to \pm \sqrt{2\mu k^*}$.
The phase portrait is shown in Figure~\ref{claudio_fig2}.

\placefigure{claudio_fig2}

\begin{figure}[!h]
\centering
\includegraphics[scale=.95]{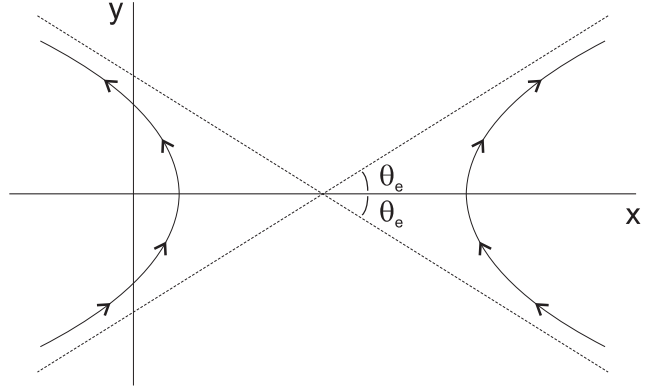}
\caption{Behavior of the solutions of the two-charged problem
with  $k^* > 0, C < 0$.}
\label{claudio_fig3}
\end{figure}

The ordinary differential equations obtained from the Hamiltonian (\ref{eq-4}) can be integrated, as in the Kepler's problem. Defining $u={\rho}^{-1}$, and using both constants of motion $k^*$ and $l$, we have
$$ d \theta =
-\frac{du}{\sqrt{\frac{2\mu k^*}{l^2}-\frac{2\mu \vert C \vert
u}{l^2}-u^2}},
$$
thus
$$
\theta - \theta' = -\arccos \left ( \frac{\frac{2 \mu \vert C
\vert}{l^2}+2u} {\sqrt{\frac{4\mu^2C^2}{l^2}+\frac{8\mu
k^*}{l^2}}} \right ),
\label{eq-theta}
$$
where $\theta'$ is a constant of integration. Therefore, $\rho$ and $\theta$ are related through the
expression
\begin{equation}
\frac{1}{\rho} = c\ [-1+e \cos (\theta - \theta')],
\label{eq-6}
\end{equation}
where
$$
c = \frac{\mu \vert C \vert }{l^2}, \enskip e =
\sqrt{1+\frac{2l^2k^*}{\mu C^2}}.
$$
The equation (\ref{eq-6}) defines a hyperbola of eccentricity $e$, whose focus is located 
at the origin. The hyperbola is sketched in Figure~\ref{claudio_fig3}; by simplicity we have chosen $\theta' = 0$. Since $e>1$ we have 
that $1>\frac{1}{e}= \cos (\theta_{e})>0$ for some $\theta_{e} \in (0,\frac{\pi}{2})$, thus the asymptotes of the hyperbola are associated to the angles $\theta_e$, $2\pi - \theta_e$.

\placefigure{claudio_fig3}

\section{The restricted charged three-body problem}

\quad In this Section we introduce the restricted charged
three-body problem, as a limiting case of the general charged
three-body problem. This restricted problem will appear when the mass of the third body tends to zero.

We remember that the Hamiltonian corresponding to the system
of bodies $i=1,2,3$ of charge $q_i$ and mass $m_i$, in an inertial frame of reference, that interact under the gravitational and Coulomb
forces on the plane, is
$$
\begin{array}{l}
\displaystyle H({\bf r}_1, {\bf r}_2, {\bf r}_3, {\bf p}_1, {\bf p}_2, {\bf
p}_3) = \\ [1pc]
\displaystyle  \sum_{i=1}^{3} \frac{\|{\bf p}_i\|^2}{2m_i} - \sum_{1 \le i < j \le 3} \frac{\lambda_{ij}}{\|{\bf r}_i- {\bf r}_j\|},
\end{array}
\label{ham-3-cargas}
$$
where
$$
\lambda_{ij}= G m_i m_j- k q_i q_j .
\label{eq-lambdas}
$$
As before, $k$ is the Coulomb constant, $G$ is the universal
gravitational constant, and ${\bf r}_i \in \mathbb{R}^2$, ${\bf p}_i \in
\mathbb{R}^2$, $i=1,2,3$ are the positions and momentum of the bodies
$1$, $2$ and $3$, respectively.

We introduce the real parameters $\alpha_i$ and $C_{ij}$ by means of
$$
\begin{array}{c}
q_i = \alpha_i {m_i}, \quad i=1, 2, 3, \\ [1pc]
C_{ij} = G-k\alpha_i\alpha_j, \quad i,j=1, 2, 3,
\end{array}
\label{def-alphas} 
$$
with the restriction $\alpha_3 \ne 0$. Notice that this last condition guarantees the existence of charge
of the third particle, therefore a different problem from the classic one.

The differential equations associated to the problem of three charged bodies are given by
\begin{equation}
\begin{array}{l}
\displaystyle \ddot{{\bf r}}_1 = - \frac{\lambda_{12}({\bf r}_1- {\bf
r}_2)}{m_1 \|{\bf r}_1 - {\bf r}_2 \|^{3}}
 - \frac{\lambda_{13}({\bf r}_1- {\bf r}_3)}{m_1 \|{\bf r}_1- {\bf
 r}_3\|^{3}},\\[1pc]
\displaystyle \ddot{\bf r}_2 = - \frac{\lambda_{21}({\bf r}_2- {\bf r}_1)}{m_2
\|{\bf r}_2- {\bf r}_1\|^{3}} - \frac{\lambda_{23}({\bf
r}_2- {\bf r}_3)}{m_2 \|{\bf r}_2- {\bf r}_3\|^{3}},\\[1pc]
\displaystyle \ddot{\bf r}_3 = - \frac{\lambda_{13}({\bf r}_3- {\bf r}_1)}{m_3
\|{\bf r}_3- {\bf r}_1\|^{3}} - \frac{\lambda_{23}({\bf r}_3- {\bf
r}_2)}{m_3 \|{\bf r}_3- {\bf r}_2\|^{3}}.
\label{eq-7}
\end{array}
\end{equation}
Taking $m_3 \rightarrow 0$ we have that $q_3
\rightarrow 0$, then $\lambda_{13} \rightarrow 0$ and
$\lambda_{23} \rightarrow 0$, therefore, by (\ref{eq-7}),
the position vectors ${\bf r}_1$ and ${\bf r}_2$ obey to
one problem of two charged bodies, that is
$$
\begin{array}{l}
\displaystyle \ddot{\bf r}_1 = - \frac{C_{12} m_2 ({\bf r}_1- {\bf
r}_2)}{\|{\bf
r}_1- {\bf r}_2\|^{3}},\\[1pc]
\displaystyle \ddot{\bf r}_2 = - \frac{C_{12} m_1 ({\bf r}_2- {\bf
r}_1)}{\|{\bf r}_1- {\bf r}_2\|^{3}}.
\end{array}
$$
On the other hand, for the third position vector we have
$$
\ddot{\bf r}_3 = - \frac{C_{13} m_1({\bf r}_3- {\bf r}_1)}{\|{\bf
r}_3- {\bf r}_1\|^{3}} - \frac{C_{23} m_2({\bf r}_3- {\bf
r}_2)}{\|{\bf r}_3- {\bf r}_2\|^{3}}.
$$
The dynamics of the third particle depends on the
motion previously chosen for the bodies $1$ and $2$, and on
the parameters
$$
C_{13}= G- k \alpha_1 \alpha_3, \quad C_{23}= G-k \alpha_2 \alpha_3.
$$
As we saw in the previous Section, in
order to have circular orbits for the primaries, the parameters must satisfy
\begin{equation}
 C_{12} = G- k \alpha_1 \alpha_2 > 0.
\label{eq-8}
\end{equation}

{\noindent \bf Remark 3.1}
{\it
Without the restriction (\ref{eq-8}) on
the parameters $\alpha_1$ and $\alpha_2$, or $m_1$, $m_2$,
$q_1$ and $q_2$, the primaries cannot follow circular trajectories, thus the circular restricted
charged three-body problem will not be well defined.}

We take a circular solution for the bodies $1$ and $2$ (inertial frame), where
both particles move around their center of mass with constant angular frequency
$\omega$. We set ${\bf r}_1 (t)= c_1 (-\cos \omega
t, -\sin \omega t) $ as solution for the first body, and
${\bf r}_2(t) = c_2 (\cos \omega t, \sin
\omega t)$ for the second one; $c_1$ and $c_2$ are
positive constants that satisfy the condition $m_1c_1=m_2c_2$, which corresponds to define the origin of the
system at the center of mass of the bodies $1$ and $2$. Since the bodies 
$1$ and $2$ have a circular trajectory, it is
appropriate to take a rotating system whose frequency is $\omega$.
The position vectors in the inertial and rotating frames are related
by means of the rotation
$$
R(t) = \left( {\begin{array}{*{20}c}
   {\cos \omega t } & {- \sin \omega t} \\ [6pt]
   {\sin \omega t} & {\cos \omega t}
\end{array}} \right)
$$
as follows:
$$
{\bf r}_1(t)= R(t) {\bf s}_1,\quad {\bf r}_2(t)= R(t) {\bf
s}_2,\quad {\bf r}_3= R(t) {\bf r}.
$$
It is not difficult to see that ${\bf s}_1$ and ${\bf s}_2$ lie on the
horizontal axis (rotating frame), that is ${\bf s}_1 = c_1
(-1,0)$, ${\bf s}_2 = c_2 (1,0) $. On the other hand, the equation that describes the movement of
the third body in the synodic frame is given by
\begin{equation}
\ddot{\bf r} = \omega^{2} {\bf r} + 2 \omega J \dot{\bf r} -
\frac{C_{13}  m_1  ({\bf r}- {\bf s}_1)}{\|{\bf r} -
 {\bf s}_1\|^{3}} - \frac{C_{23} m_2  ({\bf r}- {\bf s}_2)}{\|{\bf r}- {\bf s}_2\|^{3}},
 \label{eq-10}
\end{equation}
with
$$
J= 
\left(
\begin{array}{cc}
0 & 1\\
-1 & 0
\end{array}
\right).
$$
The second order differential equation (\ref{eq-10}) can be written
in terms of a pair of first order vector differential equations in
$\mathbb{R}^2$, that is
$$
\begin{array}{l}
\displaystyle \dot{\bf r} = {\bf v},\\[1pc]
\displaystyle \dot{\bf v} = 2 \omega J {\bf v} + \frac{\partial
\Omega}{\partial {\bf r}},
\label{sistema-restringido}
\end{array}
$$
where the function $\Omega$ is given by
\begin{equation}
\Omega = \frac{1}{2} \omega^2 {\bf r}^2 + \frac{C_{13}
m_1}{\rho_1} + \frac{C_{23} m_2}{\rho_2},
\label{eq-11}
\end{equation}
with $\rho_1 = \|{\bf r} - {\bf s}_1\|$ and $\rho_2 =
 \|{\bf r} - {\bf s}_2\|$.

\begin{figure}[!h]
\centering
\includegraphics[scale=0.8]{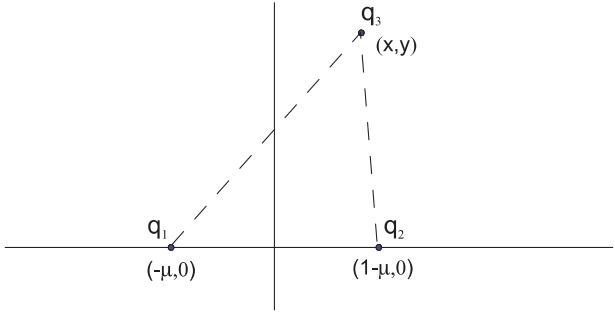}
\caption{The restricted circular charged three-body problem in
a rotating frame with convenient units.}
\label{claudio_fig4}
\end{figure}

It is convenient to introduce the parameters
$$ \mu_1 = G m_1, \quad \mu_2 =  G m_2, \quad \alpha_j \sqrt{\frac{k}{G}}=
\widetilde{\alpha_j}, \quad j=1,2,3,
$$
together with
$$
\delta = \textrm{sgn}(\alpha_3), \quad \beta_j = 1 - \widetilde{\alpha_j} \widetilde{\alpha_3}, \quad j=1,2.
\label{def-betas}
$$
Notice that four independent physical magnitudes exist:
time, distance, mass and charge. We choose $\omega$ as
unitary frequency (time), the separation between $1$ and
$2$ as unitary distance. Also units of charge and mass such that
$$
\vert \widetilde{ \alpha_3} \vert = 1,
\quad \mu_1+\mu_2 = 1.
$$
According to the selection of units of mass, and by the symmetry of the problem, it is enough to consider only one mass 
parameter belonging to $(0,1/2]$. For this we use $\mu$, where $0< \mu \leq 1/2$, $\mu_2 = \mu$ and $\mu_1 = 1-\mu$. With this we have
$$
\begin{array}{c}
C_{13}m_1 = (1-\mu)(1-\widetilde{\alpha_1}\delta) = (1- \mu) \beta_1, \\[10pt]
C_{23}m_2 = \mu(1-\widetilde{\alpha_2}\delta) = \mu \beta_2.
\end{array}
\label{eq-C13-C23}
$$
Using the condition $(1-\mu)c_1 = \mu c_2$ (center of mass), the relation $c_1+c_2=1$ (units of 
distance), and $G m_1 = 1 - \mu$, $Gm_2 = \mu$ (units of mass), it is obtained 
that ${\bf s}_1=(-\mu, 0)$ and ${\bf s}_2=(1-\mu, 0)$. The positions of the three bodies, in the synodic frame, are shown in Figure~\ref{claudio_fig4}.

\placefigure{claudio_fig4}

\begin{figure}[!h]
\centering
\includegraphics[scale=0.75]{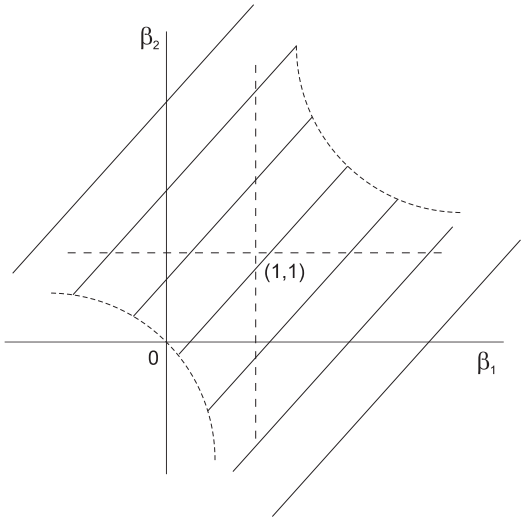}
\caption{The hatched part, delimited by 
$\beta_2=\frac{1}{\beta_1-1} + 1$, defines the allowed values of the
parameters $\beta_1$ and $\beta_2$.}
\label{claudio_fig5}
\end{figure}

The parameters that appear in (\ref{eq-11}) have been written in
terms of the new parameters $\beta_1$, $\beta_2$, but still
without the restriction associated to $C_{12}$. For this units, 
the condition $C_ {12} > 0$ is equivalent to
$\widetilde{\alpha_1}\widetilde{\alpha_2}<1$, which implies
\begin{equation}
(\beta_1-1)(\beta_2-1) < 1.
\label{eq-12}
\end{equation}
The region defined by (\ref{eq-12}) is shown in Figure \ref{claudio_fig5}.

\placefigure{claudio_fig5}

{\noindent \bf Remark 3.2}
{\it
The difference between the circular restricted charged three-body problem
and the circular restricted three-body problem of Celestial
Mechanics (which is given by $(\beta_1, \beta_2) = (1,1) $) is in
the values that can assume $\beta_1$ and $\beta_2$.
}

In the following, it is described the relation between the forces that the particle $i=1,2$,
exerts on the third particle, according to the values
that $\beta_i$ can assume:

\begin{itemize}
\item $\beta_i < 0$. The Coulomb force is repulsive and
predominates over the gravitational one.

\item $\beta_i = 0$. The gravitational and Coulomb forces have equal magnitude, being the
second one repulsive.

\item $0 < \beta_i < 1$. The Coulomb force is repulsive and the
gravitational force predominates.

\item $\beta_i = 1$. The Coulomb force is null.

\item $1 < \beta_i$. The Coulomb force is attractive.
\end{itemize}

Denoting the coordinates of ${\bf r} $ by $(x, y)$, and taking
$\omega=1$, we have that (\ref{eq-10}) is
equivalent to
\begin{equation}
\begin{array}{l}
\ddot x - 2 \dot y - x= V_x, \\[1pc]
\ddot y + 2 \dot x - y= V_y,
\end{array}
 \label{eq-13}
\end{equation}
where
\begin{equation}
V = \frac{\beta_1 (1-\mu)}{\rho_1} + \frac{\beta_2
\mu}{\rho_2},
\label{potencial-V}
\end{equation}
with
$$\rho_1= \sqrt{(x+ \mu)^2+ y^2}, \quad \rho_2= \sqrt{(x-1+\mu)^2+ y^2}.
$$
Considering the coordinates
$$
p_x= \dot x- y,\quad p_y= \dot y+  x,
$$
we obtain that
(\ref{eq-13}) is equivalent to the system
\begin{equation}
\begin{array}{l}
\displaystyle \dot x= y+ p_x = H_{p_x}, \quad \dot p_x=  V_x+ p_y= -H_x ,\\[1pc]
\displaystyle \dot y= -x+ p_y= H_{p_y}, \quad \dot p_y=  V_y- p_x= -H_y,
\label{eq-14}
\end{array}
\end{equation}
whose Hamiltonian function is
\begin{equation}
H(x, y, p_x, p_y)= \frac{1}{2}(p_x^2+ p_y^2)+ (y p_x- x p_y)- V.
\label{funcion-Hamiltoniana}
\end{equation}
The Hamiltonian function (\ref{funcion-Hamiltoniana}), the potential function (\ref{potencial-V}), and the restriction (\ref{eq-12}), define the restricted charged three-body problem. \\

{\noindent \bf Remark 3.3}
{\it
\cite{dionysiou-1} considered the circular restricted
three-charged-body problem with a different parameterization from
the one used in this work. On the other hand, \cite{rad-1,rad-2},
\cite{kunitsyn-1}, \cite{lukyanov-1},
\cite{simmons}, studied the restricted photogravitational problem 
making use of a parameterization similar to the one we have used in our work. In the photogravitational case it
is considered, in addition to the gravitational force, a repulsive (or null) force 
of magnitude inversely proportional to the square of the
distance, equivalent to the restriction $\beta_i \in (-\infty,
1]$.}

A relevant consequence of (\ref{eq-12}) is that the equilibrium solutions are confined 
to stay in a certain region of the configuration space, as we shall see in the following Section.

\section{Equilibrium solutions}

\begin{figure}[!h]
\centering
\includegraphics[scale=0.85]{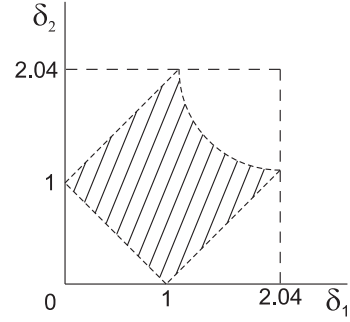}
\caption{Allowed region, on the parameters space $(\delta_1,\delta_2)=(\beta_1^{1/3},\beta_2^{1/3})$, for the existence of triangular equilibrium solutions. The region is delimited by the functions $\delta_2=\delta_1+1$, $\delta_2=\delta_1-1$, $\delta_2=1-\delta_1$, 
$\displaystyle \delta_2^3 = \frac{1}{\delta_1^3-1}+1$.}
\label{claudio_fig6}
\end{figure}

\begin{figure}[!h]
\centering
\includegraphics[scale=0.85]{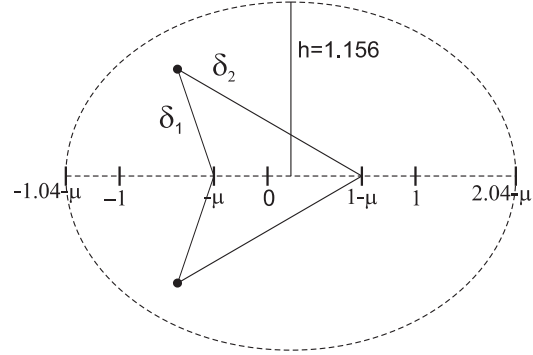}
\caption{Allowed region, on the configuration space, for the existence of triangular equilibrium solutions.}
\label{claudio_fig7}
\end{figure}

According to (\ref{eq-14}), the equilibrium solutions of the circular restricted
charged three-body problem must satisfy
$$
V_x + x = 0, \quad  V_y + y = 0,
\label{eq-equilibrios}
$$
which implies
$$
\begin{array}{l}
x f - \mu (1-\mu)\left(\frac{\beta_1}{\rho_1^3} -
\frac{\beta_2}{\rho_2^3}\right)= 0,\\[1pc]
y f = 0,
\label{eq-equilibrios-2}
\end{array}
$$
where
$$
f = 1 - \frac{\beta_1(1- \mu)}{\rho_1^3} - \frac{\beta_2\mu}{\rho_2^3}.
$$
We can separate the solutions in two types: $y \ne 0$ (non
collinear or triangular) and $y=0$ (collinear). Considering the first case, $y \ne 0$, it is required that $f=0$,
which implies that
\begin{equation}
1 - \frac{\beta_1(1- \mu)}{\rho_1^3} - \frac{\beta_2\mu}{\rho_2^3}
= 0 \label{eq-15}
\end{equation}
and
\begin{equation}
\frac{\beta_1}{\rho_1^3} - \frac{\beta_2}{\rho_2^3}=0
\label{eq-16}
\end{equation}
hold. The equation (\ref{eq-16}) requires that $\beta_1$ and $\beta_2$
have the same sign. If this happens, the sign must
be positive, so that (\ref{eq-15}) has solution. The corresponding solution 
is $\rho_1 = \beta_1^{1/3}$, $\rho_2 = \beta_2^{1/3}$. Therefore, it 
is required that the $\beta$-parameters be positive, and satisfy
the triangular inequalities $\rho_1 + 1>\rho_2$, $\rho_2 + 1>\rho_1$, $\rho_1 + \rho_2>1$. Besides this, the 
restriction (\ref{eq-12}) must be fulfilled. In order to describe the triangular solutions on the 
parameters space, we introduce $(\delta_1,\delta_2) = (\beta_1^{1/3},\beta_2^{1/3})$; we will 
use them interchangeably. In Figures \ref{claudio_fig6} and \ref{claudio_fig7} it is shown the allowed 
region for the existence of triangular equilibrium solutions, in the parameters and configuration space, respectively.

\placefigure{claudio_fig6}

\placefigure{claudio_fig7}

There are two equilibrium solutions of this type, one with $y>0$,
and another with $y<0$. It is not difficult to see that the
coordinates of these solutions are $(x,y)=(x_L,\pm y_L)$, where
$$
\begin{array}{c}
x_L = -\mu + \frac{1}{2}(\beta_1^{2/3} - \beta_2^{2/3} + 1), \\[1pc]
y_L = \frac{1}{2} \sqrt{2(\beta_1^{2/3} + \beta_2^{2/3})- (\beta_1^{2/3} - \beta_2^{2/3})^2 - 1}.
\end{array}
$$
Thus, we have proved that:\\

{\noindent \bf Theorem 4.1}~\label{teo-triangular}
{\it
There are two triangular equilibrium solutions in the circular restricted charged three-body problem, whenever $\beta_1$
and $\beta_2$ belong to the region shown in Figure~\ref{claudio_fig6}. They are given by
$$
L_4 = (x_L,y_L), \quad L_5 = (x_L,-y_L).
$$
}

{\bf Remark 4.1}
{\it
This result was proved in the corresponding models by \cite{kunitsyn-3},
\cite{lukyanov-1,lukyanov-2} and \cite{simmons}.
}\\

\begin{figure}[!h]
\centering
\includegraphics[scale=0.70]{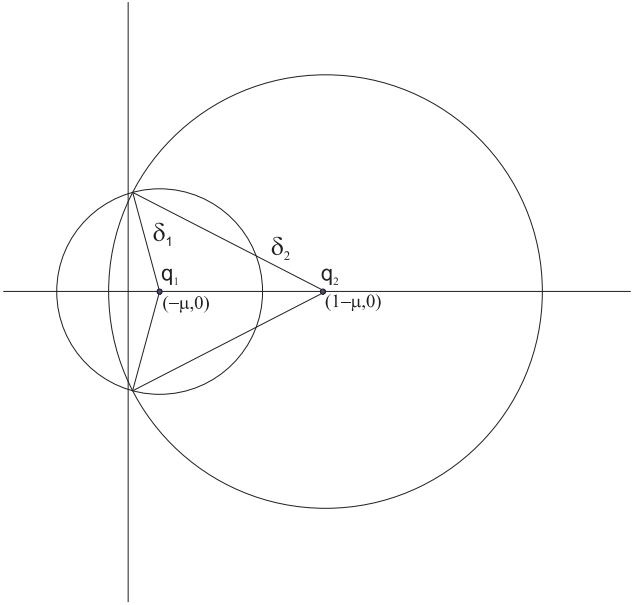}
\caption{Equilibrium solution situated to the left of the body $1$.}
\label{claudio_fig8}
\end{figure}

\begin{figure}[!h]
\centering
\includegraphics[scale=0.70]{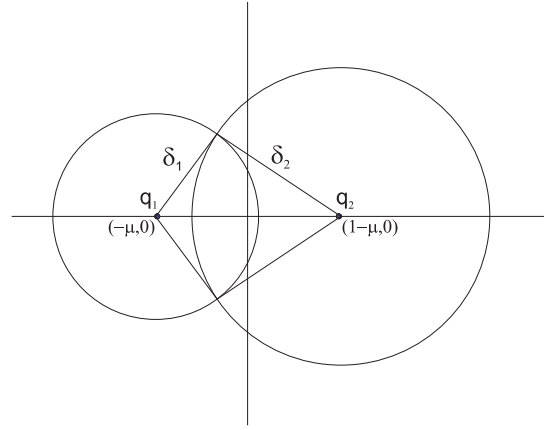}
\caption{Equilibrium solution situated between the bodies $1$ and $2$.}
\label{claudio_fig9}
\end{figure}

\begin{figure}[!h]
\centering
\includegraphics[scale=0.70]{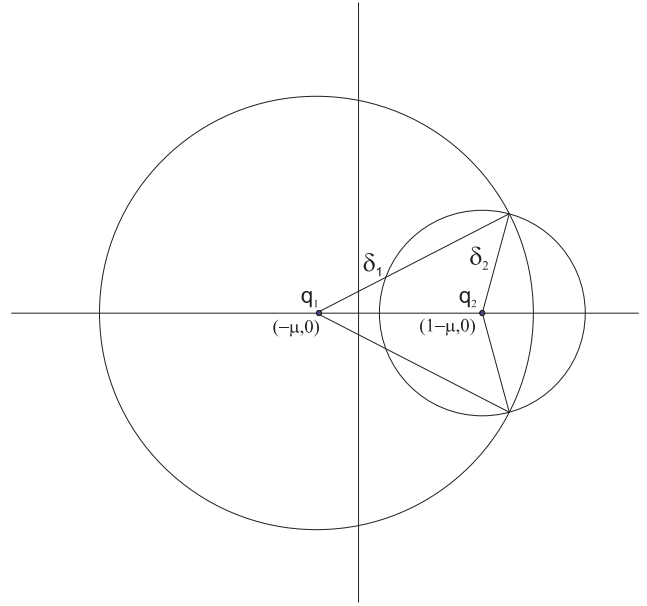}
\caption{Equilibrium solution situated to the right of the body $2$.}
\label{claudio_fig10}
\end{figure}

{\bf Proposition 4.1}~\label{prop-triangular}
{\it
The parameters $\beta_1$ and $\beta_2$ set the location of the triangular equilibrium solutions, in the following way:
\begin{itemize}
\item $\beta_2^{1/3} > \sqrt{1+\beta_1^{2/3}}$. The equilibrium
solution is located to the left of the body $1$ (see Figure~\ref{claudio_fig8}).
\item $\beta_2^{1/3} = \sqrt{1+\beta_1^{2/3}}$. The equilibrium solution is above or below the body $1$.
\item $\beta_1^{1/3} < \sqrt{1+\beta_2^{2/3}}$ or $\beta_2^{1/3}
< \sqrt{1+\beta_1^{2/3}}$. The equilibrium solution is located
between the bodies $1$ and $2$ (see Figure~\ref{claudio_fig9}).
\item $\beta_1^{1/3} = \sqrt{1+\beta_2^{2/3}}$. The equilibrium solution is above or below the body $2$.
\item $\beta_1^{1/3} > \sqrt{1+\beta_2^{2/3}}$. The equilibrium
solution is located to the right of the body $2$ (see Figure~\ref{claudio_fig10}).
\end{itemize}
}

\placefigure{claudio_fig8}
\placefigure{claudio_fig9}
\placefigure{claudio_fig10}

On the other hand, for the collinear case, i.e. $y = 0$, the relative distances
assume the form
$$\rho_1 = |x+ \mu|,\quad \rho_2 = |x+\mu-1|,$$
and the roots of the equation
\begin{equation}
F(x,\beta_1,\beta_2) = x f - \mu (1-\mu) \left(\frac{\beta_1}{\rho_1^3} -
\frac{\beta_2}{\rho_2^3}\right)=0
\label{eq-17}
\end{equation}
define equilibrium solutions (by simplicity, we also use $F(x)$ to denote $F(x,\beta_1,\beta_2)$). The problem 
of the existence of the roots of $F$ depends on $\beta_1$ and $\beta_2$, and has been solved in parametric form by
\cite{lukyanov-1}, without the condition $(\beta_1-1)(\beta_2-1)<1$. In the following we show necessary and sufficient 
conditions for the existence of the roots of $F$ under such restriction. As first step, notice 
that (\ref{eq-17}) can be rewritten as
$$
F(x) = x - \frac{\beta_1\ (1-\mu)\, (x+\mu)}{\rho_1^3} - \frac{\beta_2\ \mu\, (x+\mu-1)}{\rho_2^3}.
\label{equi-colineal-1}
$$
In order to facilitate the study of the function $F$, it is
convenient to separate the domain in three intervals:
$$
I_1 = (-\infty,-\mu), \enskip I_2 = (-\mu,1-\mu), \enskip I_3 = (1-\mu,\infty).
$$
We have the following reductions:
\begin{itemize}
\item $x \in I_1$. The relations
$x + \mu = -\rho_1$, $x + \mu - 1 = -\rho_2$ are fulfilled,
therefore $F(x)$ is reduced to $x +
\frac{\beta_1(1-\mu)}{\rho_1^2} + \frac{\beta_2\ \mu}{\rho_2^2}$.
\item $x \in I_2$. In this case
$x + \mu =
 \rho_1$, $x + \mu - 1 = -\rho_2$, then $F(x)$ becomes
 $x - \frac{\beta_1(1-\mu)}{\rho_1^2} + \frac{\beta_2\ \mu}{\rho_2^2}$.
\item $x \in I_3$. It is
verified that $x + \mu = \rho_1$, $x + \mu - 1 = \rho_2$, 
therefore $F(x)$ assumes the form $x -
\frac{\beta_1(1-\mu)}{\rho_1^2} - \frac{\beta_2\ \mu}{\rho_2^2}$.
\end{itemize}
Note that in the three intervals it is fulfilled
$$
F'(x)= 1+ 2 \beta_1 (1-\mu) \frac{1}{\rho_1^3}+ 2 \beta_2 \mu
\frac{1}{\rho_2^3},
$$
where the prime denotes differentiation respect to $x$. Next, we 
divide the plane $(\beta_1,\beta_2)$ into four quadrants and two lines, that is
$$
\begin{array}{l}
R_1 = \{ (\beta_1,\beta_2) \in \mathbb{R}^2 \enskip \vert \enskip \beta_1 > 0, \enskip \beta_2 > 0 \}, \\
R_2 = \{ (\beta_1,\beta_2) \in \mathbb{R}^2 \enskip \vert \enskip \beta_1 < 0, \enskip \beta_2 > 0 \}, \\
R_3 = \{ (\beta_1,\beta_2) \in \mathbb{R}^2 \enskip \vert \enskip \beta_1 < 0, \enskip \beta_2 < 0 \}, \\
R_4 = \{ (\beta_1,\beta_2) \in \mathbb{R}^2 \enskip \vert \enskip \beta_1 > 0, \enskip \beta_2 < 0 \}, \\
R_5 = \{ (\beta_1,\beta_2) \in \mathbb{R}^2 \enskip \vert \enskip \beta_1 = 0 \}, \\
R_6 = \{ (\beta_1,\beta_2) \in \mathbb{R}^2 \enskip \vert \enskip \beta_2 = 0 \},
\end{array}
$$
and define the $S$-regions by means of
$$
\begin{array}{l}
S_{1,1} = \{ (\beta_1,\beta_2) \in \mathbb{R}^2 \enskip \vert \enskip 0 < \beta_1 \le  1, \enskip \beta_2 > 0 \}, \\
S_{1,2} = \{ (\beta_1,\beta_2) \in \mathbb{R}^2 \enskip \vert \enskip \beta_1 >1, \enskip 0< \beta_2 < \frac{\beta_1}{\beta_1-1} \}, \\
S_2 = \{ (\beta_1,\beta_2) \in \mathbb{R}^2 \enskip \vert \enskip \beta_1 < 0, \enskip \beta_2 >  \frac{\beta_1}{\beta_1-1} \}, \\
S_{4,1} = \{ (\beta_1,\beta_2) \in \mathbb{R}^2 \enskip \vert \enskip 0 < \beta_1 < 1, \, \, \frac{\beta_1}{\beta_1-1} < \beta_2 < 0 \}, \\
S_{4,2} = \{ (\beta_1,\beta_2) \in \mathbb{R}^2 \enskip \vert \enskip \beta_1 \ge 1, \enskip \beta_2 <0 \}, \\
S_{5} = \{ (\beta_1,\beta_2) \in \mathbb{R}^2 \enskip \vert \enskip \beta_1 = 0, \enskip \beta_2 > 0 \}, \\
S_{6} = \{ (\beta_1,\beta_2) \in \mathbb{R}^2 \enskip \vert \enskip \beta_1 > 0, \enskip \beta_2 = 0 \}.
\end{array}
$$
With this, we introduce the sets $R'_i$, $i=1,\cdots,6$, which correspond respectively to $R_i$, $i=1,\cdots,6$, restricted to $(\beta_1-1)(\beta_2-1)<1$. Thus, we have
$$
\begin{array}{l}
R'_1 = S_{1,1} \cup S_{1,2}, \\
R'_2 = S_2, \\
R'_3 = \varnothing, \\
R'_4 = S_{4,1} \cup S_{4,2}, \\
R'_5 = S_{5}, \\
R'_6 = S_{6}.
\end{array}
$$
These sets define the allowed regions for $(\beta_1,\beta_2)$. Depending on the region (with exception of $R'_3$), and the interval, could exist one or two collinear equilibrium solutions. 

In the following we introduce a Theorem that relates different 
regions of the $\beta$-parameters for which the collinear equilibrium solutions exist. It shall be 
useful for demonstration of subsequent results.

{\bf Theorem 4.2}~\label{theorem-property}
{\it
Suppose that there exist collinear equilibrium solutions for $(\beta_1,\beta_2)$ defined by the inequalities
\begin{equation}
\begin{array}{c}
f_1(x_*,\mu) \le \beta_1 \le f_2(x_*,\mu), \\ [10pt]
g_1(x_*,\mu) \le \beta_2 \le g_2(x_*,\mu), \\ [10pt]
h_1(\mu) \le x_* \le h_2(\mu),
\end{array}
\label{eq-18}
\end{equation}
that is, for each pair $(\beta_1,\beta_2)$ there exist $x$ such that $F(x,\beta_1,\beta_2) = 0$.
Then, there exist collinear equilibrium solutions for those $(\beta_1,\beta_2)$ which satisfy
\begin{equation}
\begin{array}{c}
g_1(-x_*,1-\mu) \le \beta_1 \le g_2(-x_*,1-\mu), \\ [10pt]
f_1(-x_*,1-\mu) \le \beta_2 \le f_2(-x_*,1-\mu), \\ [10pt]
-h_2(1-\mu) \le x_* \le -h_1(1-\mu).
\end{array}
\label{eq-19}
\end{equation}
}

{\bf Proof:}\\
We want to show that there exist some $x$ for the new $\beta$-parameters such that this triad defines a root of $F$.
Consider the change of variables $x = - \widetilde x$, $x_* = - \widetilde x_*$, $\mu = 1 - \widetilde \mu$, 
$\beta_1 = \widetilde \beta_2$, $\beta_2 = \widetilde \beta_1$. With this, the inequalities in (\ref{eq-18}) take the form
\begin{equation}
\begin{array}{c}
f_1(- \widetilde x_*,1- \widetilde \mu) \le \widetilde \beta_2 \le f_2(- \widetilde x_*,1- \widetilde \mu), \\ [10pt]
g_1(- \widetilde x_*,1- \widetilde \mu) \le \widetilde \beta_1 \le g_2(- \widetilde x_*,1- \widetilde \mu), \\ [10pt]
h_1(1- \widetilde \mu) \le -x_* \le h_2( 1- \widetilde \mu).
\end{array}
\label{eq-20}
\end{equation}
In the following, we also use the tilde notation for the variables related with $\widetilde x$, $\widetilde x_*$, 
$\widetilde \mu$, $\widetilde \beta_1$, $\widetilde \beta_2$. Notice that the relative distances become
$\rho_1^2 = \widetilde \rho_2^2$, $\rho_2^2 = \widetilde \rho_1^2$, therefore $F(x,\beta_1,\beta_2)=0$ 
can be written as $F(\widetilde x, -\widetilde \beta_1, -\widetilde \beta_2)=0$. The later 
equation and (\ref{eq-20}) describe another set of $\beta$-parameters for which collinear equilibrium
solutions exist. Removing the tilde in these relations we obtain (\ref{eq-19}).

\begin{figure}[!h]
\centering
\includegraphics[scale=1]{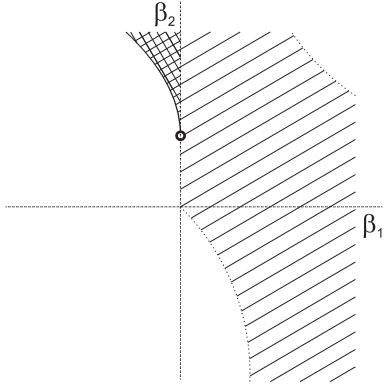}
\caption{Values of the $\beta$-parameters for which there exist collinear equilibrium 
solutions, on the interval $I_1$. Each pair $(\beta_1,\beta_2)$ belonging to the single hatched part allows only 
one equilibrium solution. On the other hand, the points associated to the double hatched part allow two different 
equilibrium solutions. Moreover, there exists only one equilibrium solution for 
the $\beta$-parameters belonging to the continuous border of the double hatched region; these parameters
are associated to a root of $F$ with multiplicity 2.}
\label{claudio_fig11}
\end{figure}

\begin{figure}[!h]
\centering
\includegraphics[scale=1]{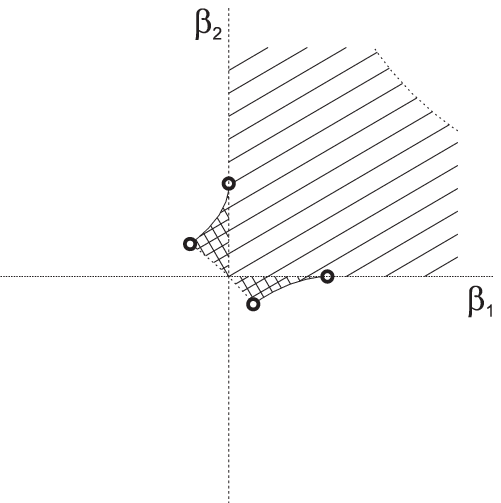}
\caption{Same description of Figure \ref{claudio_fig11}, on the interval $I_2$.}
\label{claudio_fig12}
\end{figure}

\begin{figure}[!h]
\centering
\includegraphics[scale=1]{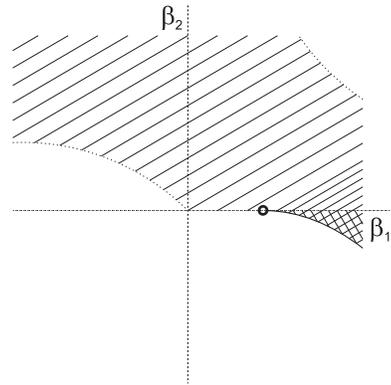}
\caption{Same description of Figure \ref{claudio_fig11}, on the interval $I_3$.}
\label{claudio_fig13}
\end{figure}

Now we give two Theorems about the necessary and sufficient conditions for the existence of the 
collinear equilibrium solutions, in terms of the $\beta$-parameters and $I_i$, $i=1,2,3$. The 
first Theorem deals with one simple root of $F$, whereas the 
second one is related to two roots. The regions of existence of the collinear equilibrium 
solutions are shown in Figures \ref{claudio_fig11}, \ref{claudio_fig12} and \ref{claudio_fig13}. \\

{\bf Theorem 4.3}~\label{teo-simple}
{\it
\begin{enumerate}
\item Region $R'_1$. There exists exactly one collinear equilibrium solution 
for $x \in I_1$, $I_2$, $I_3$.

\item Region $R'_2$. There exists exactly one collinear equilibrium solution 
for $x \in I_3$.

\item Region $R'_4$. There exists exactly one collinear equilibrium solution 
for $x \in I_1$.

\item Region $R'_5$. There exists exactly one collinear equilibrium solution 
for $x \in I_3$. On the other hand, if $\beta_2 > 1$ holds, then there exists exactly one collinear equilibrium solution for $x \in I_1$. In a similar way, if $0 <\beta_2 < 1$ is satisfied, then there exists exactly one collinear equilibrium solution for $x \in I_2$.

\item Region $R'_6$. There exists exactly one collinear equilibrium solution 
for $x \in I_1$. On the other hand, if $\beta_1 > 1$ holds, then there exists exactly one collinear equilibrium solution for $x \in I_3$. In a similar way, if $0 <\beta_1 < 1$ is satisfied, then there exists exactly one collinear equilibrium solution for $x \in I_2$.

\end{enumerate}
}

{\bf Proof:}\\
According to Theorem 4.2, items 2 and 4 imply 3 and 5 respectively, therefore it is enough to demonstrate items 1, 2 and 4.

First item, interval $I_1$. We notice that if $x \to -\infty$ then $F (x) \to -\infty$, and if $x \to -\mu$ then $F (x) \to \infty$, because $\beta_1 >0$. It is clear that the function is strictly increasing since $F'(x)>0$. The conclusion follows by the continuity of the function. A similar argument can be used to demonstrate the same item, intervals $I_2$ and $I_3$, and fourth item, interval $I_3$.

Fourth item, interval $I_1$. In this case $F(x) = x + \frac{\beta_2 \mu}{\rho_2^2}$, and the function takes its maximum value at $x=-\mu$. Also notice that if $x \to - \infty$ then $F(x) \to - \infty$, and that $F(x)$ is strictly increasing. Therefore, the function could have only one root. For the existence of the root it is required $F(-\mu) > 0$, which implies $\beta_2 > 1$. A similar argument can be used to demonstrate the same item, interval $I_2$.

Second item. We write the function $F(x)$ in terms of $\rho_2 = x + \mu - 1$, as a quotient of polynomials defined for $\rho_2 \in (0, \infty)$. The equilibrium solutions are defined by the positive roots of the fifth-degree polynomial in the numerator:
$$
\begin{array}{c}
\rho_2^5 + \rho_2^4(3-\mu) + \rho_2^3(3-2\mu) + \\ [10pt]
 + \rho_2^2((1-\mu)(1-\beta_1)-\beta_2 \mu) + \\ [10pt]
- 2\beta_2 \mu \rho_2 - \beta_2 \mu = 0.
\end{array}
$$
In this polynomial, the coefficients of the quintic, quartic and cubic terms are positives, whereas the sign of 
the coefficient of the quadratic term depends on the values of $\beta_1$ and $\beta_2$, and the remain coefficients are negative. This implies that there is only one variation in the sign of the sequence of 
the coefficients, independently of the coefficient of the quadratic term, therefore we have only one positive root.

In order to state the Theorem concerning two roots of $F$ (possibly one root of multiplicity 2), we 
introduce a new parameter $x_*$; it will be used to characterize the frontier of the region of existence 
of the collinear solutions, on the parameters space (we proceed as \cite{lukyanov-1}). Once done 
that, we introduce
\begin{equation}
\begin{array}{l}
\displaystyle \beta_1^*(x_*,\mu) = \frac{(3x_*+\mu-1)(x_*+\mu)^3}{2(1-\mu)}, \\ [10pt]
\displaystyle \beta_2^*(x_*,\mu) = \frac{(3x_*+\mu)(x_*+\mu-1)^3}{2\mu}, \\
\end{array}{}
\label{eq-G}
\end{equation}
$G(x_*,\mu) = \beta_1^*(x_*,\mu)\beta_2^*(x_*,\mu) - \beta_1^*(x_*,\mu) - \beta_2^*(x_*,\mu)$,\vspace{.3cm}
\noindent 
and $x_{r1}(\mu)$, $x_{r2}(\mu)$, to be defined in an implicit way. For fixed $\mu$, the functions $x_{r1}(\mu)$, $x_{r2}(\mu)$ 
are specific roots of the eighth-degree polynomial $G(x_*,\mu)$ in $x_*$, for instance $G(x_{r1}(\mu),\mu) = 0$. The 
second root is defined, in terms of the first one, by the expression $x_{r2}(\mu) = -x_{r1}(1-\mu)$. These functions satisfy
$-\mu <x_{r1}(\mu)<-\frac{\mu}{3}$, $\frac{1}{3}(1-\mu) <x_{r2}(\mu)<1-\mu$. At the Appendix we 
give an approximation of the roots $x_{r1}(\mu)$ and $x_{r2}(\mu)$, using regular perturbation theory. \\

The functions $x_{r1}(\mu)$ as the value of $x_*$ where the 
curve $\beta_1 = \beta_1^*(x_*)$, $\beta_2 = \beta_2^*(x_*)$, $-\mu <x_*<-\frac{\mu}{3}$ meets
$\beta_1\beta_2 - \beta_1 - \beta_2 =0$, that is $G(x_{r1}(\mu),\mu)=0$, we have $x_* \in (-\mu, x_{r1}(\mu))$ for the upper bound.

and $x_{r1}(\mu)$, $x_{r2}(\mu)$ as specific roots of the eighth-degree polynomial $G(x_*,\mu)$ in $x_*$. At the Appendix we give an approximation of the roots $x_{r1}(\mu)$ and $x_{r2}(\mu)$, using regular perturbation theory. \\

{\bf Theorem 4.4}~\label{teo-double}
{\it
\begin{enumerate}

\item  Region $R'_2$. There are at most two collinear equilibrium solutions in each one of the intervals 
$I_1$, $I_2$. A necessary and sufficient condition for the existence of the equilibrium collinear solutions 
with $x \in I_1$ is $-\beta_1^*(x_*) \le \beta_1 <0$, $\beta_2^*(x_*) \le \beta_2$, $x_* < -\mu$. On the 
other hand, the necessary and sufficient condition for $x \in I_2$ is $\beta_1^*(x_*) \le \beta_1 <0$, 
$\displaystyle \frac{\beta_1^*(x_*)}{\beta_1^*(x_*) - 1} < \beta_2 \le \beta_2^*(x_*)$, $-\mu < x_* < x_{r1}(\mu)$.

\item  Region $R'_4$. There are at most two collinear equilibrium solutions in each one of the intervals 
$I_3$, $I_2$. A necessary and sufficient condition for the existence of the equilibrium collinear solutions 
with $x \in I_3$ is $-\beta_2^*(x_*) \le \beta_2 <0$, $\beta_1^*(x_*) \le \beta_1$, $x_* > 1-\mu$. On the 
other hand, the necessary and sufficient condition for $x \in I_2$ is $\beta_2^*(x_*) \le \beta_2 <0$, 
$\displaystyle \frac{\beta_2^*(x_*)}{\beta_2^*(x_*) - 1} < \beta_1 \le \beta_1^*(x_*)$, $x_{r2}(\mu) < x_* < 1 - \mu$.
\end{enumerate}
}

{\noindent \bf Remark 4.2}
{\it
There are two types of collinear equilibrium solutions involved in Theorem 4.4. In one 
case the root of $F$ is simple, whereas in the other one the root has multiplicity 2. The 
roots of multiplicity 2 only happen when the equalities for both $\beta_1$ and $\beta_2$ 
are fulfilled. For instance, $\beta_1 = -\beta_1^*(x_*)$, $\beta_2 = \beta_2^*(x_*)$, $x_* < -\mu$, in 
the first case of $R'_2$.}\\

{\bf Proof:}\\
First item, interval $I_1$. We write the function $F(x)$ in terms of $\rho_1 = - x - \mu$, as a quotient of 
polynomials, defined for $\rho_1 \in (0, \infty)$. The equilibrium solutions are 
defined by the positive roots of the polynomial of fifth degree in the numerator:\\
$$
\begin{array}{c}
-\rho_1^5-(2+\mu)\rho_1^4 -(1+2\mu)\rho_1^3 + \\ [10pt]
+ (\beta_1(1-\mu) + (\beta_2-1)\mu )\rho_1^2 +  \\ [10pt]
+ 2\beta_1(1-\mu)\rho_1  +\beta_1(1-\mu) = 0.
\end{array}
$$
Since $\beta_1 < 0$, possibly with exception of the coefficient of the quadratic term, the coefficients are negative, 
so it is possible to have two or none sign variations in the sequence of the coefficients, therefore we have three 
options: two simple roots, one root of multiplicity 2, or none root. This implies that $F(x)$ is a concave function, so there exists a local extrema $x_*$, that is $F'(x_*)=0$.
Notice that the condition $F(x_*) \ge 0$ guarantees the existence of at least one root. Thus, for 
a given $\mu$, we want to determine those values of $\beta_1$ and $\beta_2$ such that both $F'(x_*)=0$ 
and $F(x_*) \ge 0$ are fulfilled. With this aim, as first step we will 
introduce a change of variables $\beta_1,\beta_2 \to x_*,k$ such that $x_*$ be a root of $F'(x)$ (in this case
$x_*$ defines a local minimum). Next, we will determine the allowed values for $k$; to do that we 
use $F(x_*) \ge 0$ and  $(\beta_1,\beta_2) \in R_2$. Finally, we will take into account the 
condition $(\beta_1-1)(\beta_2-1)<1$ which gives rise to $(\beta_1,\beta_2) \in R'_2$. Let $x_*<-\mu$ the 
root of $F'(x)$. We introduce the change of variables
$$
\begin{array}{c}
\displaystyle \beta_1 = \frac{1+k}{2(1-\mu)} (x_*+\mu)^3, \\ [10pt]
\displaystyle \beta_2 = \frac{k}{2\mu} (1-x_*-\mu)^3,
\end{array}
$$
which satisfies $F'(x_*)=0$. The condition $F(x_*) \ge 0$ implies $k \ge -(3x_* + \mu)$, therefore
$\beta_1$ and $\beta_2$ must be greater or equal to certain minimum values. Notice that $k = -(3x_* + \mu)$ 
defines the lower bound of $\beta_1$ and $\beta_2$, that is $-\beta_1^*(x_*)$ and $\beta_2^*(x_*)$, respectively. Besides this, 
due to $(\beta_1,\beta_2) \in R_2$, we require $\beta_1<0$, so at this point the inequalities 
$-\beta_1^*(x_*) \le \beta_1 < 0$, $\beta_2^* (x_*) \le \beta_2$, $x_*<-\mu$ are satisfied. Since $x_*<-\mu$, we 
have $\beta_2 > 1$; using this and $\beta_1 < 0$ we conclude that the inequality $(\beta_1-1)(\beta_2-1)<1$ holds,
therefore $(\beta_1,\beta_2) \in R'_2$.

Second item, interval $I_3$. It is consequence of the first item, interval $I_1$, and Theorem 4.2

First item, interval $I_2$. We proceed as for the demonstration of the first item, interval $I_1$. We write the 
function $F(x)$ in terms of $\rho = \frac{1}{x+ \mu} - 1$, as a quotient of polynomials, defined for $\rho \in (0, \infty)$. We
focus on the positive roots of fifth-degree polynomial in the numerator: \\

$-\beta_1(1-\mu)\rho^5 - 3\beta_1(1-\mu)\rho^4 +
((\beta_2-1)\mu-3\beta_1(1-\mu))\rho^3 +
((1-\beta_1)(1-\mu)+3\beta_2 \mu)\rho^2 + 3\beta_2\mu\rho +
\beta_2\mu = 0.$\\

By hypothesis $\beta_1 < 0$ and $ \beta_2 > 0$, therefore the coefficients present positive signs, possibly with exception of the coefficients of the quadratic and cubic terms. It is possible to have two or none
sign variations in the sequence of the coefficients, therefore we have three options: two simple roots, one root of 
multiplicity 2, or none root, so $F(x)$ is a concave function. Let $x_* \in (-\mu,1-\mu)$ the root of $F'(x)$. We 
use the change of variables
$$
\begin{array}{c}
\displaystyle \beta_1 = -\frac{1+k}{2(1-\mu)} (x_*+\mu)^3, \\ [10pt]
\displaystyle \beta_2 = \frac{k}{2\mu} (1-x_*-\mu)^3,
\end{array}
$$
so $F'(x_*)=0$. From $F(x_*) \le 0$ we obtain $k \le -(3x_* + \mu)$. Since $(\beta_1,\beta_2) \in R_2 \subset R'_2$ we
require $k>0$. Using both previous inequalities we obtain $x_*<-\frac{\mu}{3}$. Taking into account $F(x_*) \le 0$, it is 
concluded that $\beta_1$ must be greater or equal to a minimum value, whereas $\beta_2$ must be lesser or equal to a
maximum value, namely $\beta_1^*(x_*)$ and $\beta_2^*(x_*)$, respectively. Due to $(\beta_1,\beta_2) \in R_2$, we require 
$\beta_1 < 0$, $\beta_2 >0$, therefore we have obtained $\beta_1^*(x_*) \le \beta_1 < 0$, $0 < \beta_2 \le \beta_2^*(x_*)$ for 
$-\mu < x_* <-\frac{\mu}{3}$. The final step is to consider $(\beta_1-1)(\beta_2-1)<1$, or in an equivalent way 
$\beta_1\beta_2 - \beta_1 - \beta_2 <0$. Notice that the $\beta$-parameters in which we are interested have 
the curve $(\beta_1^*(x_*),\beta_2^*(x_*))$, $-\mu < x_* <-\frac{\mu}{3}$, as upper bound, and $\beta_1\beta_2 - \beta_1 - \beta_2 =0$ as 
lower bound. Defining $x_{r1}(\mu)$ as the value of $x_*$ where the 
curve $\beta_1 = \beta_1^*(x_*)$, $\beta_2 = \beta_2^*(x_*)$, $-\mu <x_*<-\frac{\mu}{3}$, meets
$\beta_1\beta_2 - \beta_1 - \beta_2 =0$, that is $G(x_{r1}(\mu),\mu)=0$, we have $x_* \in (-\mu, x_{r1}(\mu))$ for the 
upper bound. With the aim of obtaining a single inequality, for the region of existence, of the equilibrium solutions, we 
parameterize the lower bound, namely $\beta_1\beta_2 - \beta_1 - \beta_2 =0$, using $x_* \in (-\mu,x_{r1}(\mu))$; 
we only must consider the part of the curve that goes from the point 
$(\beta_1^*(x_{r1}(\mu)),\beta_2^*(x_{r1}(\mu)))$ to the origin. The corresponding parameterization 
is $\beta_1 = \beta_1^*(x_*)$, $\beta_2 = \frac{\beta_1^*(x_*)}{\beta_1^*(x_*)-1}$, 
$x_* \in (-\mu,x_{r1}(\mu))$. Therefore, the region of existence of the collinear equilibrium solutions 
is defined by the inequalities
$$
\begin{array}{c}
\displaystyle \beta_1^*(x_*) \le \beta_1 < 0, \\ [10pt]
\displaystyle \frac{\beta_1^*(x_*)}{\beta_1^*(x_*)-1} < \beta_2 \le \beta_2^*(x_*),  \\ [10pt]
\displaystyle - \mu < x_* < x_{r1}(\mu).
\end{array}
$$

Second item, interval $I_2$. It is consequence of the first item, interval $I_2$, and Theorem 4.2. We identify 
$x_{r2}(\mu)$ with $-x_{r1}(1-\mu)$, therefore $-x_{r1}(1-\mu)$ must be a root of $G(x_*,\mu)$. According 
to (\ref{eq-G}), we have that $\beta_1^*(-x_*,1-\mu) = \beta_2^*(x_*,\mu)$ and 
$\beta_1^*(-x_*,1-\mu) = \beta_2^*(x_*,\mu)$, which implies $G(-x_*,1-\mu)=G(x_*,\mu)$. From this we conclude 
$G(-x_{r1}(1-\mu),1-\mu) = G(x_{r1}(\mu),\mu)=0$, as required, for consistency.

\placefigure{claudio_fig11}

\placefigure{claudio_fig12}

\placefigure{claudio_fig13}

There exist collinear equilibrium solutions which have a very
simple algebraic expression. These solutions are the limit of the
triangular ones. \\

{\bf Theorem 4.5}~\label{teo-sol-particular-colineal}
{\it
Assuming that $(\beta_1,\beta_2) \in R'_1$, we have the following:
\begin{enumerate}
\item $(-\mu- \beta_1^{1/3}, 0)$ is a collinear equilibrium solution 
on the interval
$I_1$, whenever $\beta_2^{1/3}- \beta_1^{1/3}=1$.
\item $(-\mu + \ \beta_1^{1/3}, 0)$ is a collinear equilibrium
solution on the interval $I_2$, whenever  $\beta_1^{1/3}+
\beta_2^{1/3}=1$.
\item $(-\mu + \beta_1^{1/3}, 0)$ is a
collinear equilibrium solution on the interval $I_3$,
whenever $\beta_1^{1/3}- \beta_2^{1/3}=1$.
\end{enumerate}
}

{\bf Proof:}\\
Due to the similarity of the proof of these items, we will prove only
item 2. Since $x = -\mu + \beta_1^{1/3}$ we have $\rho_1 =
\beta_1^{1/3}$. Notice that $\rho_2 = \beta_2^{1/3}$ holds, 
since $\beta_1^{1/3}+ \beta_2^{1/3}=1$. Thus, the condition
on the relative distances $\rho_1 + \rho_2 = 1$ is true. Finally,
$$
\begin{array}{l}
F(-\mu + \beta_1^{1/3})=-\mu + \beta_1^{1/3}-\beta_1^{1/3} (1-\mu)
+ \beta_2^{1/3}  \mu \\ [1pc]
= \mu(-1 +\beta_1^{1/3} + \beta_2^{1/3})=0.
\end{array}
$$

\section{Stability of the equilibrium solutions}

\begin{figure}[!h]
\centering
\includegraphics[scale=0.9]{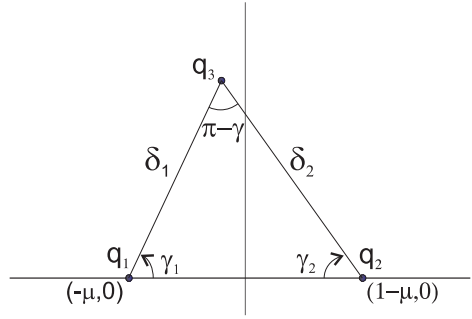}
\caption{Angles $\gamma_1$ and $\gamma_2$.}
\label{claudio_fig14}
\end{figure}

We know that the linear stability of the equilibrium solutions is
determined by the eigenvalues of the matrix
$$
A = J\ \mbox{Hess}\, H(x^*, y^*, p_x^*, p_y^*),
$$
where $(x^*, y^*, p_x^*, p_y^*)$ is an equilibrium solution of the system
(\ref{eq-14}) (see \cite{MeHa}). Then,
\begin{equation}
A = \left(
\begin{array}{cccc}
0 & 1 & 1 & 0\\
-1 & 0 & 0 & 1\\
V_{xx} & V_{xy} & 0 & 1\\
V_{xy} & V_{yy} & -1 & 0
\end{array}
\right),
\label{eq-21}
\end{equation}
where
\begin{equation}
\begin{array}{l}
\displaystyle V_{xx}=  -\frac{\beta_1 (1-\mu)}{\rho_1^3}-\frac{\beta_2 
\mu}{\rho_2^3}+ \\[1pc]
\displaystyle + 3 \beta_1 (1-\mu) \frac{(x+\mu)^2}{\rho_1^5}+ 3
\beta_2 \mu \frac{(x+\mu-1)^2}{\rho_2^5}, \\[1pc]
\displaystyle V_{xy}=  3 \beta_1 (1-\mu) \frac{(x+\mu) y}{\rho_1^5}+ 3
\beta_2 \mu \frac{(x+\mu-1) y}{\rho_2^5}, \\[1pc]
\displaystyle V_{yy}=  -\frac{\beta_1 (1-\mu)}{\rho_1^3}-\frac{\beta_2 
\mu}{\rho_2^3}+ \\[1pc]
\displaystyle + 3 \beta_1 (1-\mu) \frac{y^2}{\rho_1^5}+ 3
\beta_2 \mu \frac{y^2}{\rho_2^5}.
\end{array}
\label{eq-22}
\end{equation}
Using the relations $\rho_{1} = \beta_1^{1/3}$, $\rho_{2} =
\beta_2^{1/3}$, which are true for triangular equilibrium
solutions, or the collinear equilibrium solutions given by
Theorem 4.5, (\ref{eq-22}) is reduced to
\begin{equation}
\begin{array}{l}
\displaystyle V_{xx}=  -1 + \frac{3(1-\mu)(x + \mu )^2}{\beta_1^{2/3}} + \frac{3
\mu (x + \mu -1 )^2}{\beta_2^{2/3}}, \\[1pc]
\displaystyle V_{xy}=  \frac{3(1-\mu)(x + \mu )y}{\beta_1^{2/3}} + \frac{3 \mu
(x
+\mu -1 )y}{\beta_2^{2/3}}, \\[1pc]
\displaystyle V_{yy}=  -1 + \frac{3(1-\mu)y^2}{\beta_1^{2/3}} + \frac{3 \mu
y^2}{\beta_2^{2/3}}.
\end{array}
\label{eq-23}
\end{equation}
In order to study the characteristic polynomial associated to the matrix (\ref{eq-21}), in each equilibrium solution of
(\ref{eq-14}), it is convenient to use $\Omega = \frac{1}{2}{\bf r}^2 + V$ (see 
equations (\ref{eq-10}) and (\ref{potencial-V})). The characteristic polynomial is
\begin{equation}
\lambda^4 + \lambda^2(4-\Omega_{xx}-\Omega_{yy}) +
\Omega_{xx}\Omega_{yy}- \Omega_{xy}^2 = 0.
\label{eq-24}
\end{equation}
The triangular equilibrium solutions have been located using the
sides of the triangle: $\rho_1$, $\rho_2$ and $1$. Now we use the parameterization 
considered by \cite{lukyanov-2}, which 
consists of angles $\gamma_i$, $i=1,2$, and the unitary side; $\gamma_i$ is 
the angle between the side $\rho_i$ and $1$, as it is shown in Figure~\ref{claudio_fig14}.

\placefigure{claudio_fig14}

By means of using (\ref{eq-23}), the characteristic polynomial (\ref{eq-24}) becomes
\begin{equation}
\lambda^4 + \lambda^2 + 9\mu(1-\mu)\sin^2(\gamma)=0,
\label{eq-25}
\end{equation}
where $\gamma = \gamma_1+\gamma_2$, $\gamma \in [0,\pi]$. Notice that the characteristic equation of the 
collinear equilibrium solutions, described by Theorem 4.5, is obtained by setting $\gamma = 0, \pi$. It is straightforward 
to solve (\ref{eq-25}) using the change of variable $u = \lambda^2$. The eigenvalues of $A$ are given by
$$
\lambda = \pm \sqrt{ \frac{ -1 \pm \sqrt{F(\mu,\gamma)}}{2} },
\label{eigenvalores}
$$
where $F(\mu,\gamma) = 1-36\mu(1-\mu)\sin^2(\gamma)$. The curve $F=0$ is defined by the equation  $36 \mu\ (1-\mu)\
\sin^2(\gamma) = 1$, with $\frac{1}{2} -
\frac{\sqrt{2}}{3} \le \mu \le \frac{1}{2}$,  $\frac{1}{9} \le
\sin^2(\gamma) \le 1$. Since $\mu^* = \frac{1}{2} -
\frac{\sqrt{2}}{3}$ is the smallest value of the mass parameter
such that the function $F$ can be zero, the associated angle
must be $\gamma = \frac{\pi}{2}$. For $\mu \in (\mu^*,\frac{1}{2}]$, always there exist 
angles $\gamma =  \gamma_\mu, \pi - \gamma_\mu$, where
$ \gamma_\mu = \arcsin ( \frac {1}{6 \sqrt{\mu(1 -\mu)}} )$, in a such way that $F(\mu,\gamma)=0$ ($\gamma_\mu$ and
$\pi - \gamma_\mu$ are different from each other). In particular, for $\mu = \frac{1}{2}$ the corresponding 
angles are $\gamma_0$ and $\pi - \gamma_0$, where $\gamma_0 = \arcsin( \frac{1}{3})$.

The stability of the triangular equilibrium solutions is determined by the value of $F$; in 
the following Theorem we state four possible cases. Although the equilibrium solutions limit 
of the triangular ones (collinear equilibrium solutions discussed in Theorem 4.5) are not triangular, 
we have included them in the fourth item of following Theorem.\\

{\bf Theorem 5.1} \label{estabilidad-triangular}
{\it
The triangular equilibrium solutions, as well as the solutions limit of the triangular ones, of the circular 
restricted charged three-body problem, satisfy:
\begin{enumerate}
\item For the values of $(\mu, \gamma)$ such that $F(\mu,\gamma) < 0$, the equilibrium 
solutions are unstable in the Lyapunov sense.
\item For the values of $(\mu, \gamma)$ such that $F(\mu,\gamma)= 0$, the equilibrium 
solutions are linearly unstable.
\item For the values $(\mu, \gamma)$ such that $0 < F(\mu,\gamma)<1$, the equilibrium 
solutions are linearly stable.
\item For the values $(\mu, \gamma)$ such that $F(\mu,\gamma)=1$, the equilibrium 
solutions are linearly unstable.
\end{enumerate}
}

{\bf Proof:}\\
In order to cover all possible cases, observe that $-8 \le F(\mu,\gamma) \le 1$ holds 
for $\mu \in (0,\frac{1}{2}]$, $\gamma \in [0,\pi]$.

First item. In this case the eigenvalues of the matrix $A$ are $\lambda = \pm\sqrt{\frac{-1 \pm
i \sqrt{|F(\mu, \gamma)|}}{2}}$. The conclusion follows since the real part of
the eigenvalues is different from zero.

Second item. We observe that the eigenvalues of the matrix $A$ are $\lambda= \pm \frac{\sqrt{2}}{2}\
i$, each one with multiplicity 2. The matrix $A$ is non diagonalizable because for each $\lambda$
its eigenvector is 
$$
\left( \frac{2\lambda+V_{xy}}{D}, \frac{\lambda^2-1-V_{xx}}{D}, 
\frac{\lambda^2+1+V_{xx}+\lambda V_{xy}}{D}, 1 \right),
$$
where
$D = \lambda^3 + \lambda(1-V_{xx}) + V_{xy}$. Thus, the equilibrium solutions are linearly unstable.

Third item. The eigenvalues of the matrix $A$ are $\lambda =\pm\ i \sqrt{\frac{1 \pm \sqrt{F(\mu,\gamma)}}{2}}$. The solution is linearly stable, since all eigenvalues are pure imaginary and distinct.

Fourth item. In this case the eigenvalues of the matrix $A$ are $\lambda = 0$ with multiplicity 2, and $\lambda = \pm i$, 
so $A$ is non diagonalizable. The conclusion follows. Notice that $F(\mu,\gamma)=1$ can only happen for $\gamma=0,\pi$.

\begin{figure}[!h]
\centering
\includegraphics[scale=0.80]{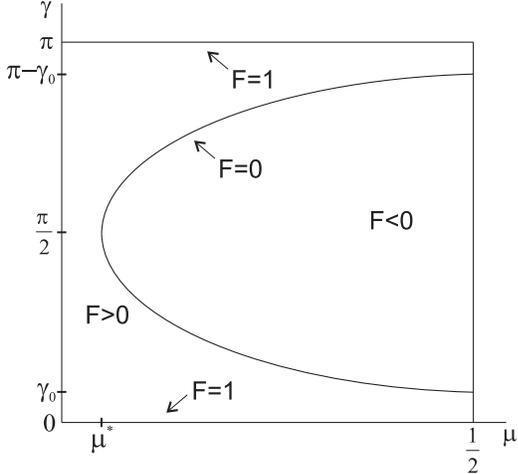}
\caption{Description of the values of $F$ as function of $(\mu,\gamma)$.}
\label{claudio_fig15}
\end{figure}

In Figure~\ref{claudio_fig15} are indicated the regions described by Theorem 5.1. Similar results were
previously obtained by \cite{lukyanov-2}.

\placefigure{claudio_fig15}

\subsection{Stable region in the configuration space}

In order to study the stable region in the configuration space, we define the restricted configuration space 
as the points $(x,y)$ which meet $\rho_1 \ne 0$, $\rho_2 \ne 0$, $\rho_1+1 \ge \rho_2$, $\rho_2+1 \ge \rho_1$, 
$\rho_1 + \rho_2 \ge 1$, $(\rho_1^3-1)(\rho_2^3-1)<1$. With this, we avoid collisions. Moreover, all the triangular
equilibrium solutions have physical sense, according to Theorem 4.1.

We want to show the evolution of the stable region, in the restricted configuration space, as $\mu$ increases. As first step, we 
apply the law of cosines to the triangle of sides $1$, $ \rho_1$, $\rho_2$, and angle of interest $\pi-\gamma$. From this we get
\begin{equation}
\rho_1^2+ \rho_2^2+2 \rho_1 \rho_2 \cos(\gamma)=1.
\label{eq-28}
\end{equation}
By geometry, we know that $\gamma =$ constant defines two arcs in the restricted configuration space, for fixed $\mu$. One arc 
satisfies $y \ge 0$, whereas the other one $y \le 0$; we will refer to them as upper and lower arcs, respectively. With the 
aim of describing the arcs, we write (\ref{eq-28}) in terms of $(x,y)$. After some algebraic manipulations, (\ref{eq-28}) becomes
\begin{equation}
\displaystyle \left(x-\frac{1}{2}+\mu\right)^2 + \left(y \pm \frac{ \cos(\gamma) }{2 \sin(\gamma)}\right)^2 = \frac{1}{4 \sin^2(\gamma)}.
\label{eq-29}
\end{equation}
The upper(lower) arc is defined by the circumference (\ref{eq-29}) with plus(minus) sign.

According to Theorem 5.1, and Figure \ref{claudio_fig15}, if $\mu \in (0,\mu^*)$ then the stable region
is characterized by $\gamma \in (0,\pi)$, and at $\mu = \mu^*$ appears the unstable region, associated to
$\gamma = \frac{\pi}{2}$, so the corresponding stable region is defined by 
$\gamma \in (0,\frac{\pi}{2}) \cup (\frac{\pi}{2},\pi)$. On the other hand, for 
fixed $\mu \in (\mu^*,\frac{1}{2}]$, the stable region is conformed by the points 
which satisfy $\gamma \in (0,\gamma_\mu) \cup (\pi-\gamma_\mu,\pi)$, where 
$\gamma_\mu = \arcsin ( \frac {1}{6 \sqrt{\mu(1 -\mu)}} )$. In the following, it is described the evolution 
of the stable region, in the restricted configuration space, according to the variation of the mass parameter.

\begin{itemize}
\item $0<\mu<\mu^*$. Exception made of the points along the horizontal axis, all the triangular equilibrium solutions
are stable, since $0 < F(\mu,\gamma) < 1$. It is outlined in Figure~\ref{claudio_fig16}.

\placefigure{claudio_fig16}

\begin{figure}[!h]
\centering
\includegraphics[scale=0.75]{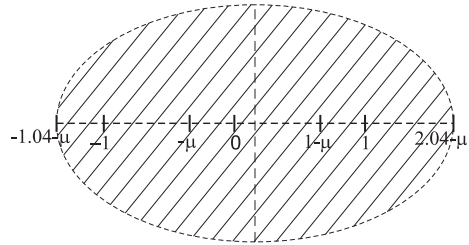}
\caption{Stable region (hatched part), for $ \mu \in (0, \mu^*)$, in the restricted configuration space.}
\label{claudio_fig16}
\end{figure}

\item $\mu = \mu^*$. The unstable region is conformed by the points on the horizontal axis, and those 
defined by $\gamma = \frac{\pi}{2}$. Except for the mentioned points, all the triangular equilibrium solutions are stable. It 
is shown in Figure~\ref{claudio_fig17}.

\placefigure{claudio_fig17}

\begin{figure}[!h]
\centering
\includegraphics[scale=0.75]{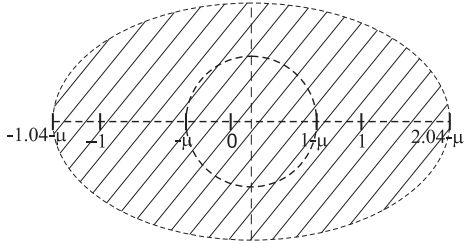}
\caption{Stable region (hatched part), for $\mu = \mu^*$, in the restricted configuration space. The dashed 
circumference indicates the unstable region associated to $\gamma = \frac{\pi}{2}$.}
\label{claudio_fig17}
\end{figure}

\item $\mu^* < \mu \le \frac{1}{2}$. The stable region, for fixed $\mu$, is defined by the points which meet
$\gamma \in (0,\gamma_\mu) \cup (\pi-\gamma_\mu,\pi)$, where $\gamma_\mu = \arcsin ( \frac {1}{6 \sqrt{\mu(1 -\mu)}} )$. It
is outlined in Figure \ref{claudio_fig18}.

\placefigure{claudio_fig18}

\begin{figure}[!h]
\centering
\includegraphics[scale=0.75]{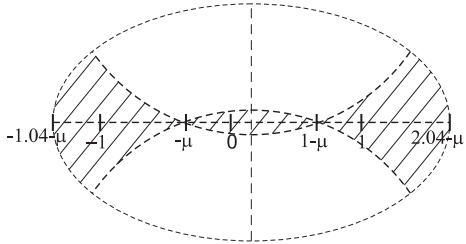}
\caption{Stable region (hatched part), for fixed $\mu \in (\mu^*,\frac{1}{2}]$, in the restricted configuration space. The dashed 
arcs are associated to $\gamma = \gamma_\mu, \pi - \gamma_\mu$.}
\label{claudio_fig18}
\end{figure}

\end{itemize}

\subsection{Stable region in the parameters space}

The evolution of the stable region also can be studied on the parameters space 
$(\delta_1, \delta_2)$. In analogy to what 
was done for the configuration space, we define the restricted parameters space as the points
$(\delta_1,\delta_2)$ which satisfy $\delta_1 \ne 0$, $\delta_2 \ne 0$, $\delta_1+1 \ge \delta_2$, $\delta_2+1 \ge \delta_1$, 
$\delta_1 + \delta_2 \ge 1$, $(\delta_1^3-1)(\delta_2^3-1)<1$. Similar results were
previously shown by \cite{simmons}, without the restriction on $\beta_1$ and $\beta_2$. 

The equation $\delta_1^2+ \delta_2^2+2 \cos(\gamma) \delta_1 \delta_2=1$ (see (\ref{eq-28})) defines an ellipse with 
center $(0,0)$, and semi-axes $\frac{1}{\sqrt{1 \pm \cos(\gamma)}}$, rotated 
$\pm \frac{\pi}{4}$ radians in a counterclockwise sense with respect to the horizontal axis (both signs hold 
due to the dependence, of the semi-axes, on the harmonic function). Notice that the ellipses
defined by $\gamma$ and $\pi-\gamma$ are related through a rotation of $\frac{\pi}{2}$ radians.

As we did for the restricted configuration space, we describe the evolution of the stable region, in the restricted parameters 
space, according to the variation of $\mu$.

\begin{itemize}
\item $0<\mu<\mu^*$. With the exception of the points that satisfy $\gamma=0,\pi$, all the triangular equilibrium solutions
are stable. It is outlined in Figure \ref{claudio_fig19}.

\placefigure{claudio_fig19}

\begin{figure}[!h]
\centering
\includegraphics[scale=0.75]{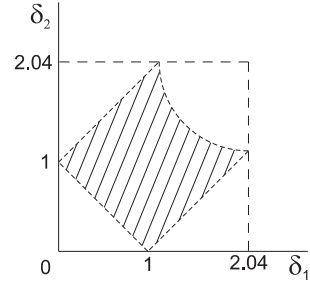}
\caption{Stable region (hatched part), for $ \mu \in (0, \mu^*)$, in the restricted parameters space.}
\label{claudio_fig19}
\end{figure}

\item $\mu = \mu^*$. The unstable region is conformed by the 
points which meet $\gamma=0,\pi,\frac{\pi}{2}$. With the exception of these points, all the triangular equilibrium solutions
are stable. It is outlined in Figure \ref{claudio_fig20}.

\placefigure{claudio_fig20}

\begin{figure}[!h]
\centering
\includegraphics[scale=0.75]{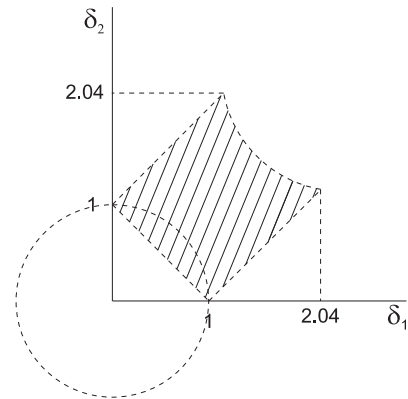}
\caption{Stable region (hatched part), for $\mu = \mu^*$, in the restricted parameters space. The dashed 
circumference indicates the unstable region associated to $\gamma = \frac{\pi}{2}$}
\label{claudio_fig20}
\end{figure}

\item $\mu^* < \mu \le \frac{1}{2}$. The stable region, for fixed $\mu$, is constituted by the points which satisfy
$\gamma \in (0,\gamma_\mu) \cup (\pi-\gamma_\mu)$, where $\gamma_\mu = \arcsin ( \frac {1}{6 \sqrt{\mu(1 -\mu)}} )$. 
It is shown in Figure \ref{claudio_fig21}.

\placefigure{claudio_fig21}

\begin{figure}[!h]
\centering
\includegraphics[scale=0.75]{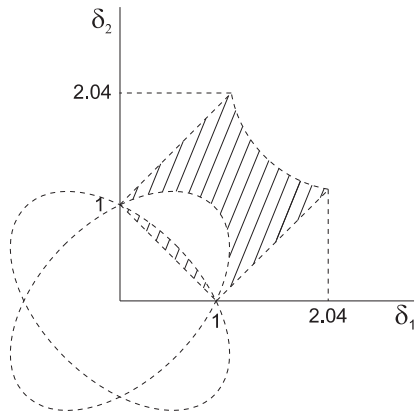}
\caption{Stable region (hatched part), for fixed $\mu \in (\mu^*,\frac{1}{2}]$, in the restricted configuration space. The curves
$\gamma = \gamma_\mu, \pi - \gamma_\mu$ are denoted with dashed ellipses.}
\label{claudio_fig21}
\end{figure}
\end{itemize}

\section{Conclusions}

The condition $(\beta_1-1)(\beta_2-1)<1$ must be fulfilled by the parameters $\beta_1$ and 
$\beta_2$, for the proper establishment of the planar circular restricted charged three-body problem. As a consequence 
of such inequality, the triangular equilibrium solutions are confined to a certain 
region of the configuration space. In a similar way, such restriction reduces the possible values of $\beta_1$ and $\beta_2$ for which 
collinear equilibrium solutions exist.

As happens in the classical restricted case, the stability of the triangular solutions is guaranteed with a small mass parameter $\mu$. It is interesting that, although the collinear equilibrium solutions (limit of the triangular ones) are linearly unstable, near of them appear stable equilibrium solutions, at 
least in the linear sense.

In this work, we have considered the triangular equilibrium solutions in the plane, as well as its linear stability, and the existence of the collinear equilibrium solutions. It would be interesting to study the stability of the equilibrium collinear solutions which are not limit of those triangular ones, as well as the equilibrium solutions in the space.

\appendix

\section{Approximations of the roots $x_{r1}(\mu)$ and $x_{r2}(\mu)$}

In order to characterize the roots $x_{r1}$, $x_{r2}$ of $G(x_*,\mu)$, we use a regular 
perturbation approach, considering $\mu$ as parameter of perturbation. Instead of work 
with $G(x_*,\mu)$, which is not well defined for $\mu = 0$, we deal with
$\widetilde G(x_*,\mu) = 4\mu(1-\mu)G(x_*,\mu)$, which is a polynomial in $\mu$. We remember that such roots satisfy 
$-\mu <x_{r1}<-\frac{1}{3} \mu$, $\frac{1}{3}(1-\mu) <x_{r2}<1-\mu$. Notice that in the limit $\mu \to 0$ we have
$\widetilde G(x_*,\mu) \to 3x_*(x_*-1)^4(2+x_*+2x_*^2+3x_*^3)$, therefore $x_{r1} \to 0$ and $x_{r2} \to 1$. For $x_{r1}$ 
we have a regular perturbation problem. We approximate $x_{r1}$ in a power series (four terms), that is
$$
x_{r1} = \sum_{i=1}^4 a_i \mu^i.
$$
Equating to zero the coefficients of the powers of $\mu$ in $\widetilde G(x_{r1},\mu)=0$, we get 
$a_1 = -\frac{1}{3}$, $a_2 = 0$, $a_3 = 0$, $a_4 = -\frac{8}{81}$, therefore the required
approximation becomes
$$
x_{r1} = -\frac{1}{3} \mu - \frac{8}{81}\mu^4.
$$
One the other hand, $x_{r2}$ cannot be handled directly by a regular perturbation approach. Nevertheless, by means of the change 
of variable $x_{r2} = 1 + \epsilon z$, where $\epsilon = \mu^{1/4}$, we obtain a regular problem with $z$ as variable, and 
$\epsilon$ as perturbation parameter. Using
$$
z = \sum_{i=0}^3 b_i \epsilon^i
$$
in $\widetilde G(1 + \epsilon z,\mu)=0$ we get $b_0 = - \sqrt[4]{\frac{4}{27}}$, $b_1 =  \frac{11}{36 \sqrt{3}}$, 
$b_2=\frac{67}{864 \sqrt[4]{12}} $, $b_3=\frac{497}{486}$, thus we obtain the approximation
$$
x_{r2} = 1 - \sqrt[4]{\frac{4}{27}} \mu^{\frac{1}{4}} + \frac{11}{36 \sqrt{3}}\mu^{\frac{1}{2}} + \frac{67}{864 \sqrt[4]{12}} \mu^{\frac{3}{4}} - \frac{497}{486}\mu.
$$
In this case, we had two possible values for $b_0$, one positive, the other one negative. In order to satisfy the inequality $x_{r2}<1-\mu$, the negative vale was chosen.

\acknowledgments The first author is pleased to acknowledge the financial support from CONACYT and PROMEP, M\'exico. The second author has been partially supported by Fondecyt 1130644.

\end{document}